\documentclass[final,3p,times]{elsarticle}
\usepackage[latin9]{inputenc}
\usepackage{array}
\usepackage{verbatim}
\usepackage{bm}
\usepackage{multirow}
\usepackage{amsmath}
\usepackage{amssymb}
\usepackage{graphicx}
\usepackage{esint}
\usepackage{amsfonts}
\usepackage{booktabs}
\usepackage{siunitx}

\usepackage{tikz}
\usetikzlibrary{shapes.geometric, arrows, positioning, matrix}

\makeatletter

\usepackage{tabularx}\usepackage{subfigure}\usepackage{subfigure}\usepackage{color}

\usepackage{setspace}

\usepackage{bm}

\journal{ }

\makeatother

\begin{document}

\title{THINC-scaling scheme that unifies VOF and level set methods}

\author[ad1]{Ronit Kumar}

\author[ad1,ad2]{Lidong Cheng}

\author[ad3]{Bin Xie}

\author[ad1]{Feng Xiao\corref{cor}}

\address[ad1]{School of Engineering, Department of Mechanical Engineering, Tokyo Institute of Technology,
\\
Tokyo, 152-8550, Japan}

\address[ad2]{School of Aerospace Engineering, Beijing Institute of Technology,\\
Beijing, 100081, China}

\address[ad3]{ School of Naval Architecture, Department of Ocean and Civil Engineering, Shanghai Jiaotong University, 
\\Shanghai, 200240, China} 

\cortext[cor]{Corresponding authors: Dr. F. Xiao (Email: xiao.f.aa@m.titech.ac.jp)}

\begin{abstract}
\indent We present a novel interface-capturing scheme, THINC-scaling, to unify the VOF (volume of fluid) and the level set methods, which have been developed as two completely different approaches widely used in various applications. The THINC-scaling scheme preserves at the samectime the advantages of both VOF and level set methods, i.e. the mass/volume conservation of the VOF method and the geometrical faithfulness of the level set method. THINC-scaling scheme allows to represent interface with high-order polynomials, and has algorithmic simplicity which eases its implementation in unstructured grids. 
\end{abstract}

\maketitle

\section{Introduction}

VOF (volume of fluid) and level set are the two most popularly used methods in capturing moving interface, and find their applications in diverse fields, such as numerical simulation of multiphase fluid dynamics, graphic processing, topological optimization and many others. Historically, they were independently developed based on completely different concepts and solution methodologies, and have their own superiority and weakness. 

VOF method\cite{hirt1981volume,youngs1982time,lafaurie1994modelling,rider1998reconstructing,scardovelli2000analytical}  uses the volume fraction of one out of multiple species in a control volume (mesh cell) to describe the distribution of the targeted fluid in space. The VOF function by definition has a value between 0 and 1, and the interfaces can be identified as isosurfaces (3D) or contours (2D) of a fractional VOF value, say 0.5 for example.  A more accurate way to represent the interface is using geometrical reconstruction, such as the PLIC  (Piecewise Linear Interface Calculation) schemes which are currently accepted as the main-stream VOF methodology. Rigorous numerical conservation can be guaranteed if a finite volume method is used to transport the VOF function, which is found to be crucial in many applications. The VOF function is usually characterized by large jump or steep gradient across the interface. So, directly using the VOF function to retrieve the geometrical information of the interface, such as normal and curvature, may result in large error. 

The level set method\cite{osher1988fronts,sethian1999level,osher2003implicit}, on other hand, defines the field function as a signed distance function (or level set function) to the interface, which possesses a uniform gradient over the whole computational domain and thus provides a perfect field function to retrieve the geometrical properties of an interface. However, the level function is not conceptually nor algorithmically conservative. The numerical solution procedure, including both transport and reinitialization,   does not guarantee the conservativeness in numerical solution. It may become a fatal problem in many applications, like multiphase flows involving bubbles or droplets.              

Efforts have been made to combine the VOF method and level set method, which lead to the coupled level set/VOF methods (CLSVOF) \cite{sussman2000coupled,menard2007coupling, yang2006adaptive, sun2010coupled}. The PLIC  type VOF method is blended with the level set method so as to improve both conservativeness and geometrical faithfulness in numerical solution. 

In this paper, we propose a new scheme that unifies the VOF and level set methods, based on the observation that the VOF field can be seen as a scaled level set field using the THINC (Tangent of Hyperbola Interface Capturing) function, which has been used in a class of schemes for capturing moving interface \cite{xiao2005simple,xiao2011revisit,ii2012interface,xie2017toward,qian2018}. The resulting scheme, so-called THINC-scaling scheme, converts the field function from level set to VOF by the THINC function, and converts the VOF field back to the corresponding level set field via an inverse THINC function. So, it can make use of the advantages of both VOF and level set at different stages of solution procedure, which eventually realizes the high-fidelity computation of moving interface regarding both numerical conservativeness and geometrical representation. 

\section{The connection between level set and VOF functions}

We consider an interface $\partial \Omega$ separating two kinds of fluids, fluid 1 and fluid 2, occupying volumes $ \Omega^1$ and $ \Omega^2$ respectively in space. We introduce the following two indicator functions to identify the different fluids and the interface. 

\begin{figure}[htbp]
	\centering
	\subfigure[] {
		\centering
		\includegraphics[width=0.4\textwidth]{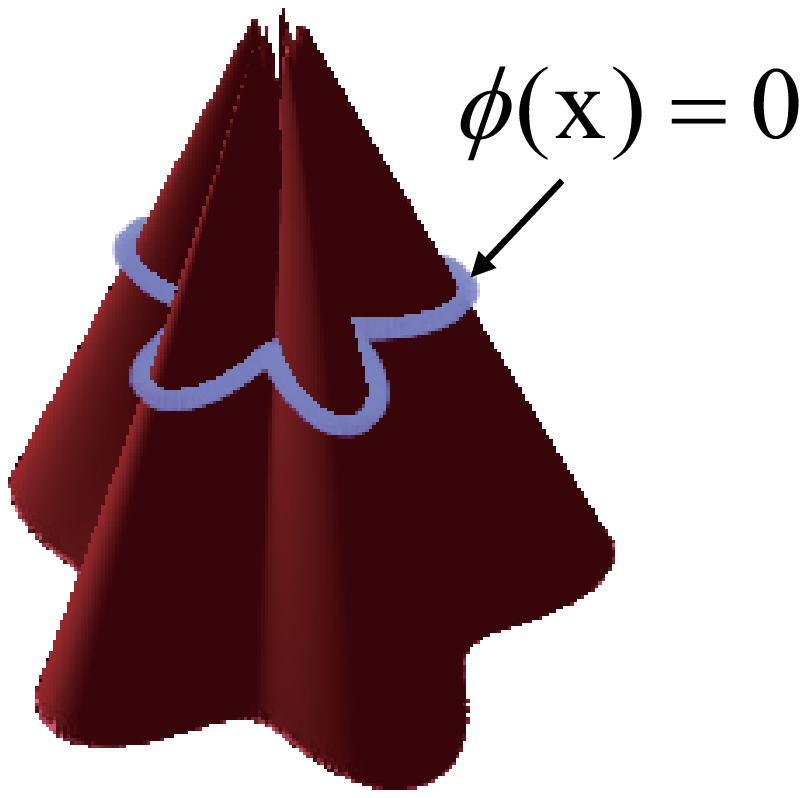}} \hspace{0.5cm}
	\subfigure[] {
		\centering
		\includegraphics[width=0.4\textwidth]{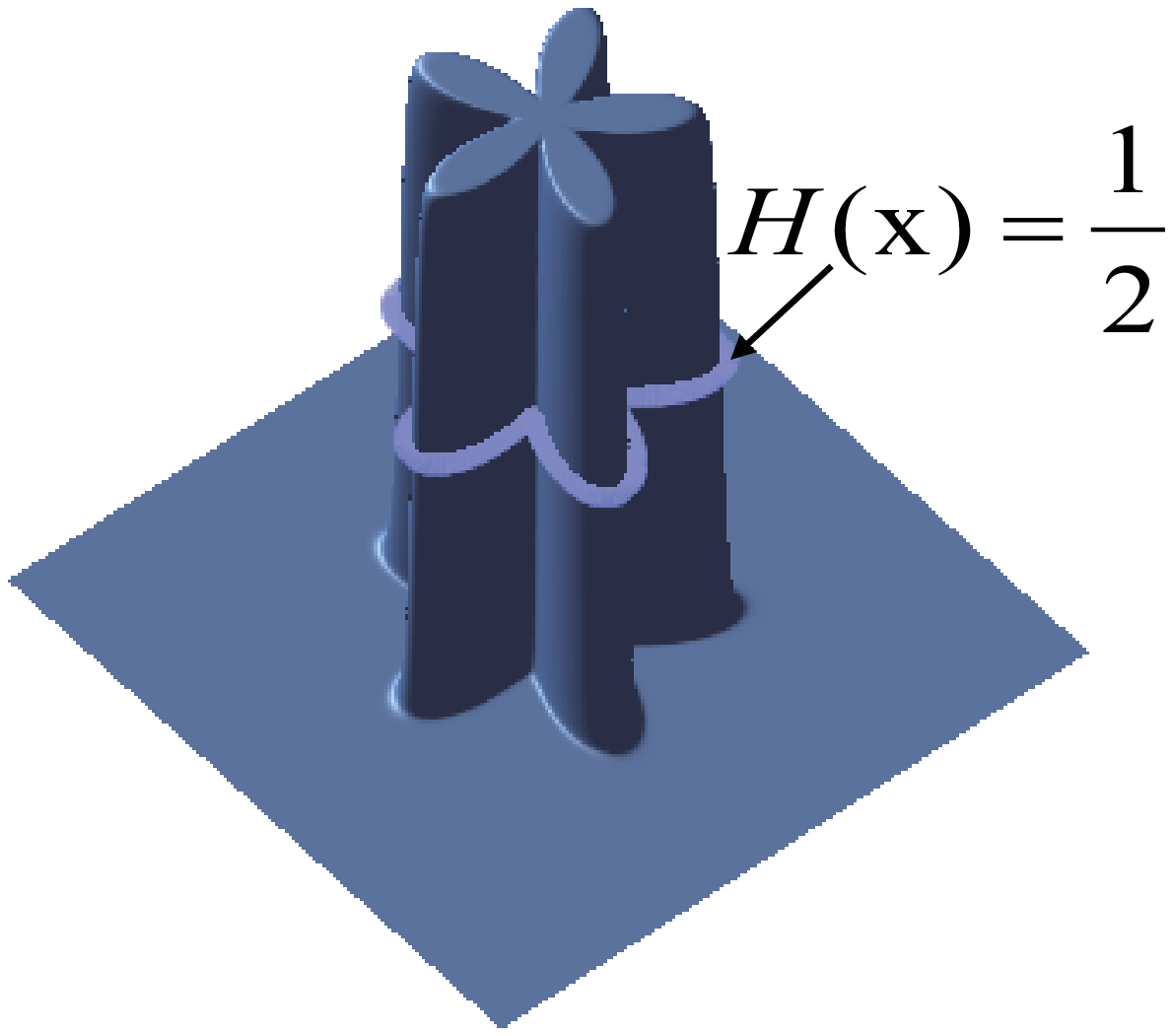}}
	\caption{Two indicator functions to identify multi-materials and interface: (a) The level set function; (b) The THINC function (a VOF function).}
\label{fun}
\end{figure}
 
\begin{itemize}
\item Level set function (Fig.\ref{fun}(a)): \\
The level set function is defined as a signed  distance from a point $\textbf{x}=(x,y,z)$ in three dimensions to the interface $\partial  \Omega$ by
\begin{equation}
\phi(\textbf{x})=\left\{
\begin{array}{rcc}
\displaystyle{\inf_{\textbf{x}_I\in \partial  \Omega}\|\textbf{x}-\textbf{x}_I\|} & & \text{if } \textbf{x}\in  \Omega^1 \\
0\quad\quad\quad & & \text{if } \textbf{x}\in\partial  \Omega \\
-\displaystyle{\inf_{\textbf{x}_I\in \partial  \Omega}\|\textbf{x}-\textbf{x}_I\|} & & \text{if } \textbf{x}\in  \Omega^2
\end{array}, \right.
\label{ls}
\end{equation}
where $\textbf{x}_I$ represents any point on the interface, and referred to as interface point. 

\item VOF function (Fig.\ref{fun}(b)): \\
The VOF function in the limit of an infinitesimal control volume is the Heaviside function. Assuming the VOF function as the abundance of fluid 1, we have the Heaviside function in its canonical form as
\begin{equation}
H^c(\textbf{x})=\left\{\begin{array}{ll}{1} & \text{if } \textbf{x}\in  \Omega^1 \\  
{\frac{1}{2}} & \text{if } \textbf{x}\in\partial  \Omega \\ 
{0} & \text{if } \textbf{x}\in  \Omega^2 \end{array}. \right.
\label{heaviside}
\end{equation}
\end{itemize}

Recall that 
\begin{equation}
H^c(\textbf{x})=\frac{1}{2} \lim_{\beta \rightarrow \infty}\left(1+\tanh \left(\beta\textbf{x})\right)\right), 
\label{hv-def}
\end{equation}
we get a continuous Heaviside function 
\begin{equation}
H(\textbf{x})=\frac{1}{2} \left(1+\tanh \left(\beta\textbf{x})\right)\right) 
\label{hv-beta}
\end{equation}
with a finite steepness parameter $\beta$. 

Given a computational mesh with cells of finite size, the cell-wise VOF function is defined by 
\begin{equation}
\bar{H}_{\Omega_l}=\frac{1}{|\Omega_l|}\int_{\Omega_l}H(\textbf{x})d\Omega,
\label{vof-def}
\end{equation}
where $\Omega_l$ is the target cell with a volume $|\Omega_l|$. In practice, we recognize the interface cell where an interface cuts through, in terms of the VOF velue, by $\varepsilon \le \bar{H}_{\Omega_l} \le 1-\varepsilon$ with $\varepsilon$ being a small positive, e.g. $\varepsilon={10}^{-8}$.

From the definition of level set function \eqref{ls}, we know that \eqref{hv-beta} scales a  level set function to a Heaviside  function.

Now, we establish a connection between the level set function and the Heaviside function. 

\begin{itemize}
\item {THINC scaling ($\phi\mapsto H)$}: 

We convert the level set field by  
\begin{equation}
H(\textbf{x},t)=\frac{1}{2} \left(1+\tanh \left(\beta\left(\mathcal {P}(\textbf{x}) \right)\right)\right), 
\label{thinc-p}
\end{equation}
where $\mathcal {P}(\textbf{x})$ is a polynomial 
\begin{equation}
\mathcal {P}(\textbf{x})=\sum_{r,s,t=0}^{p}{a}_{rst}{x}^{r}{y}^{s}{z}^{t}
\label{ls-poly}
\end{equation}
whose coefficients are determined from the level set function through the following constraints.
\begin{equation}
    \frac{\partial^D \mathcal {P}(\textbf{x})}{\partial x^{d_x}\partial y^{d_y}\partial z^{d_z}}=    \frac{\partial^D \phi (\textbf{x})}{\partial x^{d_x}\partial y^{d_y}\partial z^{d_z}}, \ \ ({d_x}, {d_y},  {d_z})=0,1,2\cdots \ \text{and } \ {d_x}+{d_y}+{d_z}=D.
    \label{ls-constraint}
    \end{equation}
In practice, we calculate the coefficients of $\mathcal {P}(\textbf{x})$  via numerical approximations using the discrete level set field available in the computational domain.   

We refer to \eqref{thinc-p} as the THINC scaling formula, and \eqref{ls-poly} as the level set polynomial.

\item {Inverse THINC scaling ($H\mapsto \phi)$}: 

Given the THINC function, we can directly compute the corresponding level set function by 
\begin{flalign} 
{\phi}(\textbf{x})=\frac{1}{\beta}{\tanh}^{-1}\left(2{H}(\textbf{x})-1\right) \  \ {\rm or} \ \ {\phi}(\textbf{x})=\frac{1}{\beta}\ln\left(\frac{{H}(\textbf{x})}{(1-{H}(\textbf{x})}\right). 
\label{THINC_inv}
\end{flalign}
Formula \eqref{THINC_inv} is referred to as the inverse THINC scaling that converts the Heaviside function to the level set function.

\end{itemize}

\begin{description}
\item {Remark 1.} Formulae \eqref{thinc-p} and \eqref{THINC_inv} provide analytical relations to uniquely convert between the level set and the THINC functions, which unifies the two under a single framework handleable with conventional mathematical analysis tool, and more importantly allows us to build interface-capturing schemes that take advantages from both VOF and level set methods. 

\item {Remark 2.} The interface is defined by  
\begin{equation}
\mathcal {P}(\textbf{x}) =0,  
\label{p0}
\end{equation}
where the level set polynomial $\mathcal {P}(\textbf{x})$ defined in \eqref{ls-poly} enables to accurately formulate the geometry of the interface. In principle, we can use arbitrarily high order surface polynomial to represent the interface straightforwardly without substantial difficulty. 

\item {Remark 3.} A finite value of the steepness parameter ${\beta}$ modifies the Heaviside function to a continuous and differentiable function \eqref{thinc-p}, which  serves an adequate approximation to the VOF function with adequately large steepness parameter $\beta$. 

\end{description}

\section{The THINC-scaling scheme for moving interface capturing}

We assume that the moving interface is transported by a velocity field $\textbf{u}$. Thus, the two indicator functions  are advected in the Eulerian form by the following equations, i.e. 
\begin{flalign}
\frac{\partial H}{\partial t}+ \nabla\cdot\left(\textbf{u}H\right)=H\nabla\cdot\textbf{u}
\label{h-adv}
\end{flalign}
for the THINC function ${H\left(\textbf{x},t\right)}$, and 
\begin{flalign}
\frac{\partial \phi}{\partial t}+ \textbf{u}\cdot\nabla\phi=0 
\label{ls-adv}
\end{flalign}
for the level set function ${\phi\left(\textbf{x},t\right)}$. 

Next, we present the THINC-scaling scheme to simultaneously solve \eqref{h-adv} and \eqref{ls-adv}. 
The computational domain is composed of non-overlapped discrete grid
cells $\Omega_{i}\,(i=1,2,\ldots,N)$ of the volume $\left|\Omega_{i}\right|$, which can be either structured or unstructured grids.  For any target cell element $\Omega_i$, we denote its mass 
center by $\textbf{x}_{ic}=(x_{ic},y_{ic},z_{ic})$, and its $J$ surface segments of areas $\left|\Gamma_{ij}\right|$ by $\Gamma_{ij}$ with $j=1,2,\ldots,J$. The outward unit normal is denoted by  ${\mathbf{n}}_{ij}=(n_{xij},n_{yij},,n_{zij})$. 

Assume that we know at time step $n$ ($t=t^n$) the VOF value   
\begin{equation}
\bar{H}_i^n=\frac{1}{|\Omega_i|}\int_{\Omega_i}H(\textbf{x},t^n)d\Omega,
\label{vofn}
\end{equation}
for each cell, and the level set value 
\begin{equation}
\phi_i^n=\phi(x_{ic},y_{ic},z_{ic},t^n)
\label{lsn}
\end{equation}
at each cell center, we use the third-order TVD Runge-Kutta scheme\cite{shu88} for time
integration to update both VOF and level set values, $\bar{H}_i^{n+1}$ and $\phi_i^{n+1}$, to the next time step $n+1$ ($t=t^{n+1}=t^{n}+\Delta t$).

We hereby summarize the solution procedure of the THINC-scaling scheme for one  Runge-Kutta sub-step that advance $\bar{H}_i^{m}$ and $\phi_i^{m}$ at sub-step $m$ to $\bar{H}_i^{m+1}$ and $\phi_i^{m+1}$ at sub-step $m+1$.   

\begin{description}
    \item {\bf Step 1.} Compute the level set polynomial of $p$th order for cell $\Omega_{i}$ from the level set field using the constraint condition \eqref{ls-constraint},
\newline
\begin{equation}
\mathcal{P}_{i}\left({\textbf x} \right)=\sum_{r,s,t=0}^{p}{a}_{rst}{X}^{r}{Y}^{s}{Z}^{t}
\end{equation}
\par where $\left(X,Y,Z\right)$ is the local coordinates with respect to the center of cell ${\Omega}_{i}$, i.e.  $X={x}-{x}_{ic}$, $Y={y}-{y}_{ic}$, $Z={z}-{z}_{ic}$. The coefficients ${a}_{rst}$ are computed by Lagrange interpolation or least square method using the level set values $\phi_i^{m}$ in the target and nearby cells. 

  \item {\textbf {Step 2}}:
Construct the cell-wise THINC function under the constraint of VOF value by 
\newline
\begin{flalign} \label{mass_constraint}
\frac{1}{|{\Omega}_{i}|}{\int}_{{\Omega}_{i}}\frac{1}{2}\left(1+\tanh\left(\beta\left(\mathcal{P}_{i}\left({\textbf x}\right)+{\phi}_{i}^{\Delta}\right)\right)\right)d \Omega={\bar H}_{i}^{m}. 
\end{flalign}
With a pre-specified $\beta$ and the surface polynomial obtained at step 1, the only unknown ${\phi}_{i}^{\Delta}$ can be computed from \eqref{mass_constraint}. In practice, we use the numerical quadrature detailed in \cite{xie2017toward}, and the resulting nonlinear algebraic function of ${\phi}_{i}^{\Delta}$ is solved by the Newton iterative method. See \cite{xie2017toward} for details. 

We then get the THINC function 
\begin{flalign} \label{thincf_mass_constraint}
H^m_{i}({\bf x})=\frac{1}{2}\left(1+\tanh\left(\beta\left(\mathcal{P}_{i}\left({\textbf x}\right)+{\phi}_{i}^{\Delta}\right)\right)\right),
\end{flalign}
which satisfies the conservation constraint of the VOF value, and ${\phi}_{i}^{\Delta}$ is a correction to the interface location due to the conservation. 
The piece-wise interface in each interface cell is defined by  
\begin{flalign} \label{cell-p}
\psi_{i}\left({\textbf x}\right)\equiv\mathcal{P}_{i}\left({\textbf x}\right)+{\phi}_{i}^{\Delta}=0. 
\end{flalign}
We refer to \eqref{cell-p} as the Polynomial Surface of the Interface (PSI) in cell ${\Omega}_{i}$. 

\item{\textbf {Step 3}}: Update the VOF function by solving \eqref{h-adv} through the following finite volume formulation, 
\begin{equation}
\bar{H}^{m+1}_{i}=\bar{H}^{m}_{i}-\frac{\Delta t}{\left|\Omega_{i}\right|}\sum_{j=1}^{J}\left(\int_{\Gamma_{ij}}\left(({\bf u} \cdot {\bf n}) H^m_{i}({\bf x})_{iup}\right)d\Gamma\right)+\frac{\bar{H}^m_{i}}{\left|\Omega_{i}\right|}\sum_{j=1}^{J}\left({\bf u} \cdot {\bf n})_{ij}\left|\Gamma_{ij}\right|\right)\Delta t,\label{fvm_semid-1}
\end{equation}
where the upwinding index $iup$ is determined by 
\begin{equation}
\begin{split}iup=\begin{cases}
=i,\ {\rm for}\ ({\bf u} \cdot {\bf n})_{ij}>0;\\
=ij,\ {\rm otherwise},
\end{cases}\end{split}
\label{eq:upwind_index}
\end{equation}
and $ij$ denotes the index of the neighboring cell that shares cell boundary 
$\Gamma_{ij}$ with target cell $\Omega_{i}$. The integration on cell surface is computed by Gaussian quadrature formula. 

\item{\textbf {Step 4}}: Update level set value at cell center using a semi-Lagrangian method as follows. 
\begin{description}
    \item {\bf Step 4.1}: 
We first find the departure point ${\textbf x}_{id}$ for each cell center ${\textbf x}_{ic}$ by solving the initial value problem, 
\begin{equation}
\begin{split}\begin{cases}
\displaystyle {\frac{{d}{\textbf x}}{d\tau}=-{\textbf u}\left({\textbf x},t^n+\tau\right)} \\
{\textbf x}(0)={\textbf x}_{ic}
\end{cases}\end{split}
\label{ivp}
\end{equation}
up to $\tau=\Delta t =t^{n+1}-t^n$, which leads to ${\textbf x}_{id}={\textbf x}(\tau)$. We use a second-order Runge-Kutta method to solve the ordinary differential equation \eqref{ivp}.
    \item {\bf Step 4.2}: Update the level set value  ${\phi}_{i}$ at cell center ${\textbf x}_{ic}$   using the Lagrangian invariant solution, 
\begin{flalign} \label{ls_sl}
{\phi}_{i}^{m+1}=\tilde{\phi}_{id}^{m}\left({\bf x}_{id}\right),   
\end{flalign}
where $\tilde{\phi}_{id}^{m}$ is the level function on cell $\Omega_{id}$ where the departure point ${\textbf x}_{id}$ falls in. 
Using the inverse THINC-scaling formula \eqref{THINC_inv}, we immediately get the level set function from \eqref{thincf_mass_constraint},
\begin{flalign} \label{ls-THINC}
\tilde{\phi}_{id}^{m}\left({\bf x}\right)=\frac{1}{\beta}{\tanh}^{-1}\left(2H^m_{id}({\bf x})-1\right), \end{flalign}
\end{description}
which gives the level set value everywhere in cell $\Omega_{id}$ that includes the departure point ${\textbf x}_{id}$. 
\item{\textbf {Step 5}}:
Reinitialize the level set field. 
We fix the level set values computed from \eqref{ls-THINC} for the interface cells which are identified by  ${\epsilon}_{1}{\leq}{\bar H}^{m}_{i}{\leq}{1-{\epsilon}_{2}}$ with ${\epsilon}_{1}$ and ${\epsilon}_{2}$ being  small positive numbers.  The level set values at the cell centers away from the interface region are reinitialized to satisfy the Eikonal equation, 
\begin{equation}
|\nabla\phi|=1. 
\end{equation}

We use the Fast Sweeping Method (FSM) \cite{zhao2005fast} on structured grid and the iteration method \cite{Dianat2017unstructured} on unstructured grid for reinitializing level set values. 

\item{\textbf {Step 6}}:
Go back to Step 1 for next sub-time step calculations.

\end{description}{}

\begin{description}
\item {Remark 4.} The THINC-scaling scheme shown above unifies the VOF method and level set method. Eq.\eqref{thincf_mass_constraint} retrieves the VOF field from the level set field, while \eqref{ls-THINC} retrieves the level set field from the VOF field with numerical conservativeness.

\item {Remark 5.} The inverse THINC-scaling \eqref{ls-THINC} facilitates a semi-Lagrangian solution without any spatial reconstruction or interpolation, such as those used in \cite{strain99}. This step essentially distinguishes the present scheme from the coupled THINC/level set method in \cite{qian2018}, where the level set function is updated by a  fifth-order Hamilton-Jacobi WENO scheme with a 3rd-order TVD Runge-Kutta time-integration scheme. 

\item {Remark 6.} The interface is cell-wisely the PSI defined by \eqref{cell-p}, $\psi_i({\textbf x})=0$,  within the interface cells, which provides the sub-cell interface structure with geometrical information, such as position, normal direction and curvature to facilitate the computation of so-called sharp-interface formulation. It distinguishes the present scheme with superiority from any other algebraic interface-capturing methods. 

\item {Remark 7.} As discussed in \cite{xiao2011revisit}, the steepness parameter ${\beta}$ can be estimated by the thickness of the jump transition across the interface using 
\begin{equation}
 \beta=\frac{1}{\eta}{\tanh}^{-1}\left(1-2\varepsilon\right)   
\end{equation}
where $\eta$ denotes the normalized half thickness of the jump with respect to the cell size, and $\varepsilon$ is a small positive number to define the range of interface transition layer in terms of the VOF value, i.e. $\varepsilon \le H(\textbf{x})\le 1-\varepsilon$, we use $\varepsilon={10}^{-8}$ in this work. In order to keep a 3-cell thickness for the interface, $\eta$ can be set as $1.5$, which  approximately results in $\beta\approx6$. Our numerical experiments show that the THINC method can resolve sharply interfaces within one or two cells by using larger $\beta$.

\end{description}

\section{Numerical tests}

We verify the THINC-scaling scheme to capture moving interfaces using some advection benchmark tests. We focus on the numerical tests in two dimensions on both structured (Cartesian) grid and unstructured (triangular) grid. 
The numerical errors in terms of the VOF field are quantified via the $L_1$ error norm $\left(\ref{L1}\right)$. 

\begin{equation}\label{L1}
E(L_1)=\sum_{i}|\bar{H}_{i}-\bar{H}_{i}^{e}||\Omega_{i}|,
\end{equation}
where $\bar{H}_{i}$ and $\bar{H}^e_{i}$ are respectively the numerical and exact VOF values.

\subsection{Solid body rotation test}
In this test, so-called Zalesak's slotted disk test \cite{zalesak1979fully}, initially a circle with a radius of 0.5 centered at $\left(0.5,0.75\right)$ in a unit square computational domain is notched with a slot defined by $\left(|x-0.5|\leq0.025 \: {\text {and}} \: y\leq0.85\right)$. The slotted circle is rotated with the velocity field given by $\left(y-0.5,0.5-x\right)$. 

The steepness parameter is set $\beta=6$, and a quadratic polynomial 
\begin{equation}
\mathcal{P}_{i}\left({\textbf x} \right)=\sum_{r,s=0}^{2}{a}_{rs}{X}^{r}{Y}^{s}
\end{equation}
is used in this test. 

\begin{figure}[htbp]
	\centering
	\subfigure[] {
		\centering
		\includegraphics[width=0.3\textwidth]{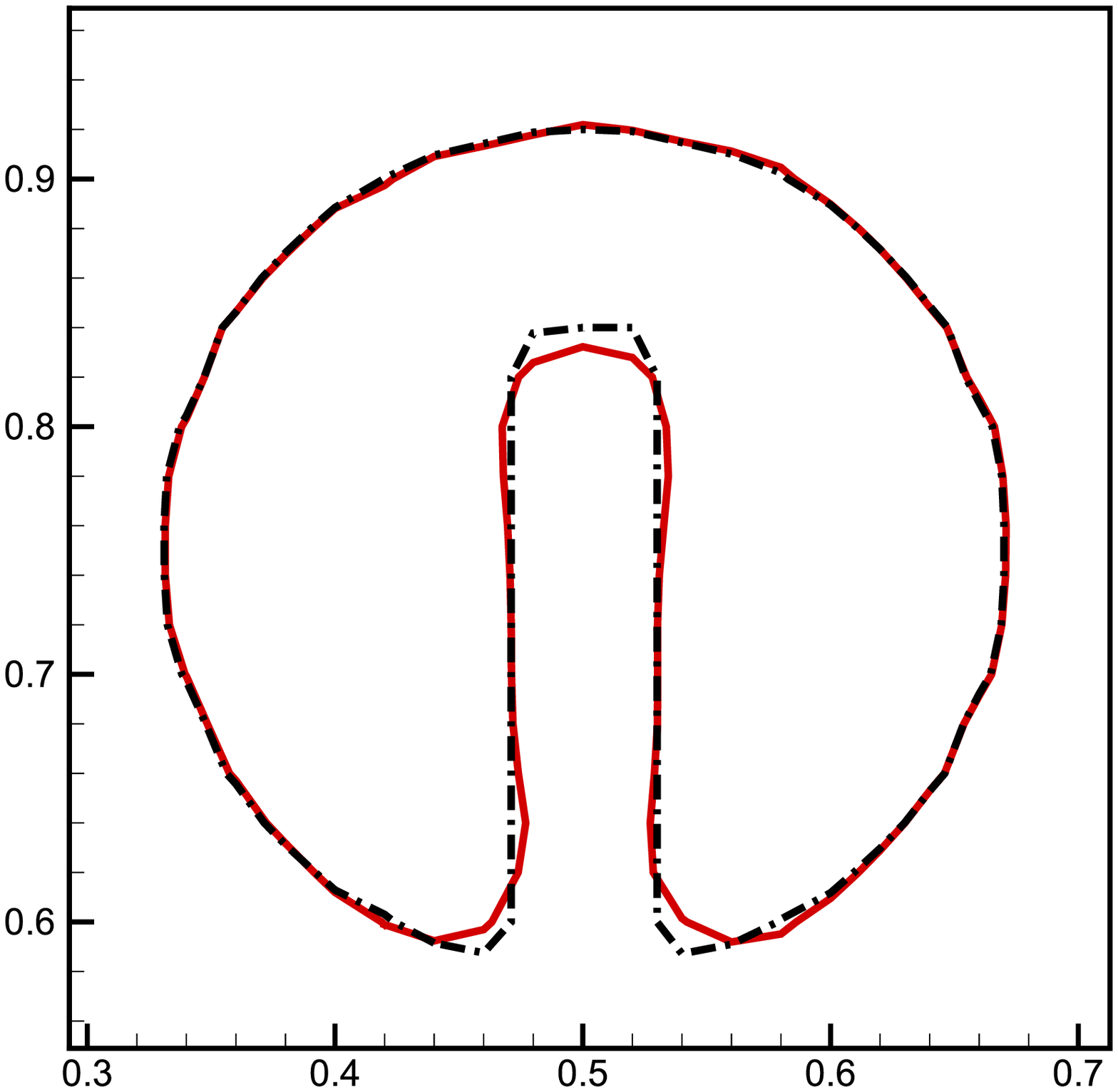}} \hspace{0.4cm}
	\subfigure[] {
		\centering
		\includegraphics[width=0.3\textwidth]{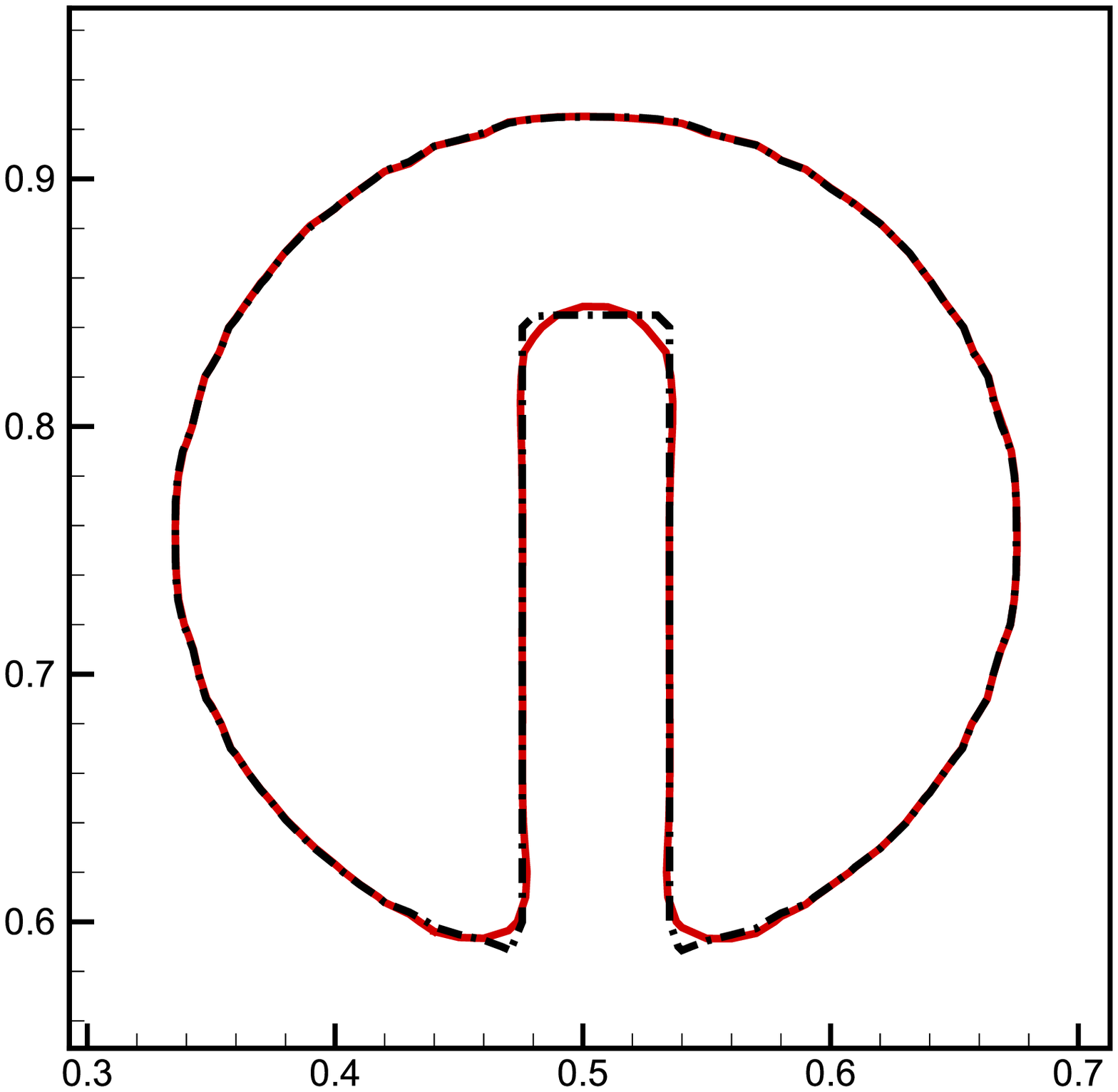}} \hspace{0.4cm}
	\subfigure[] {
		\centering
		\includegraphics[width=0.3\textwidth]{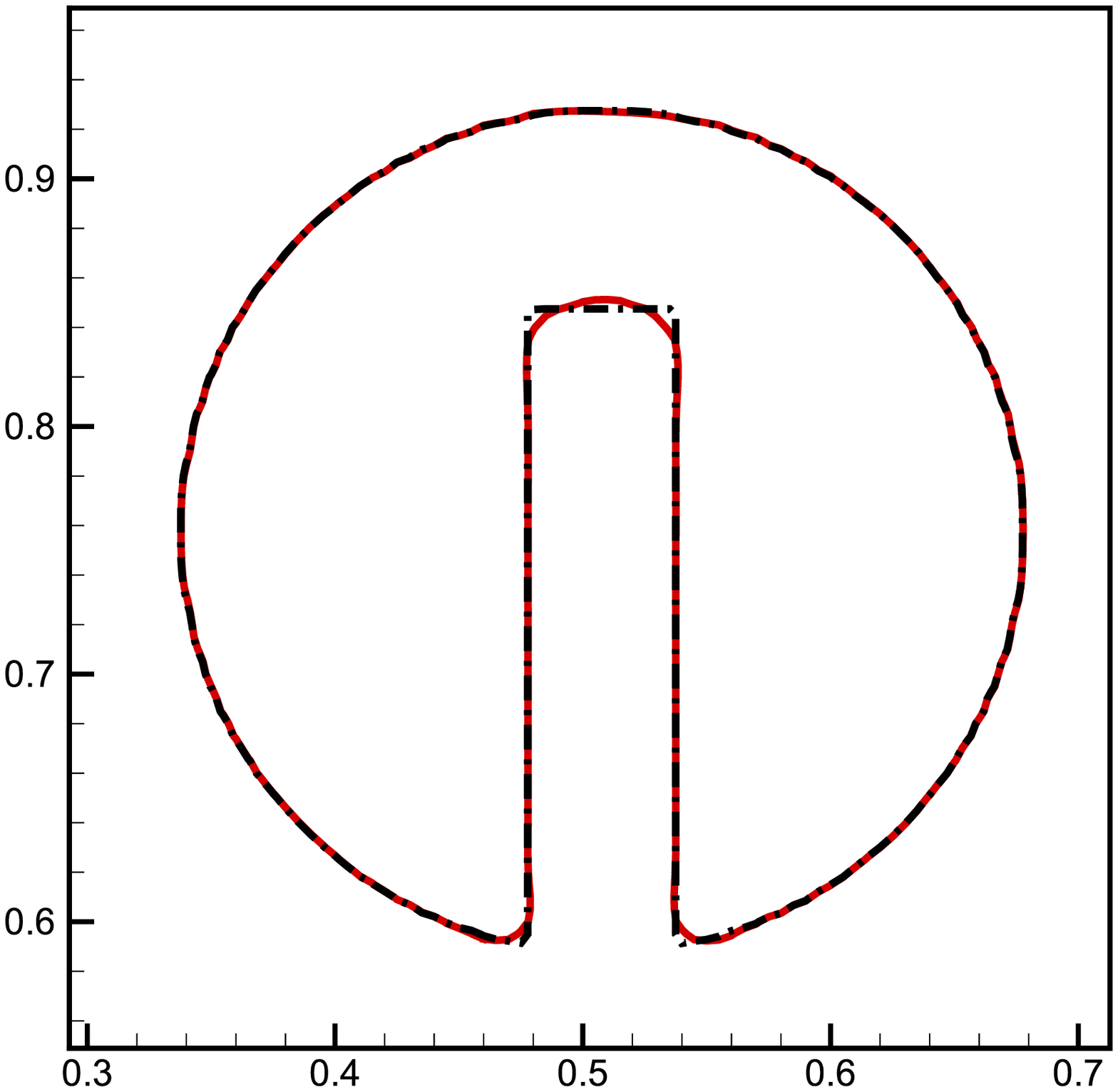}}
	\caption{VOF 0.5 contour line in Zalesak solid body rotation test after one revolution on (a) $50\times50$, (b) $100\times100$ and (c) $200\times200$ meshes. The black dashed line stands for the exact solution, and the red solid line for the numerical solution. }
\label{zalesak_results}
\end{figure}

\begin{figure}[htbp]
	\centering
	\subfigure[] {
		\centering
		\includegraphics[width=0.3\textwidth]{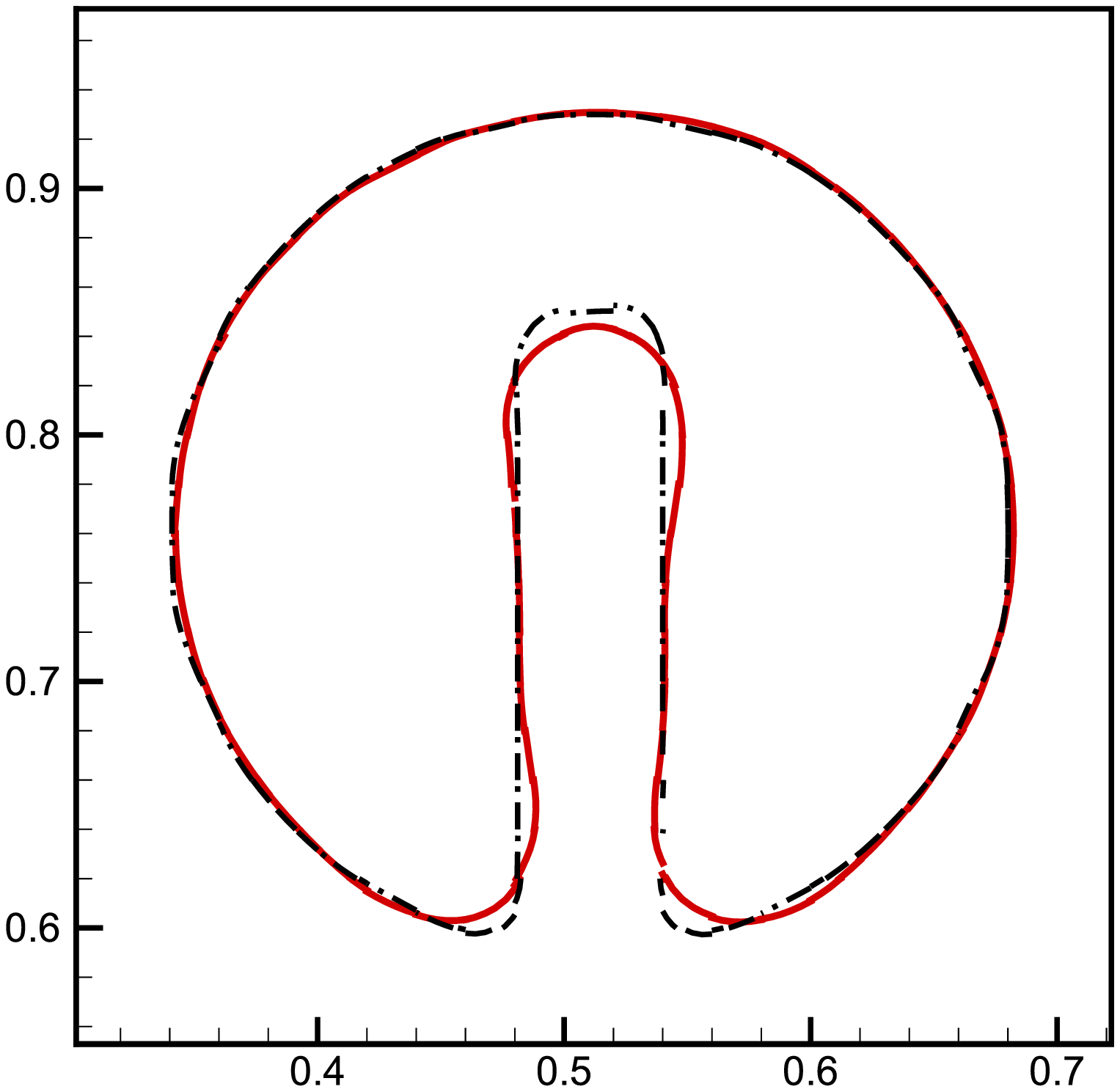}} \hspace{0.4cm}
	\subfigure[] {
		\centering
		\includegraphics[width=0.3\textwidth]{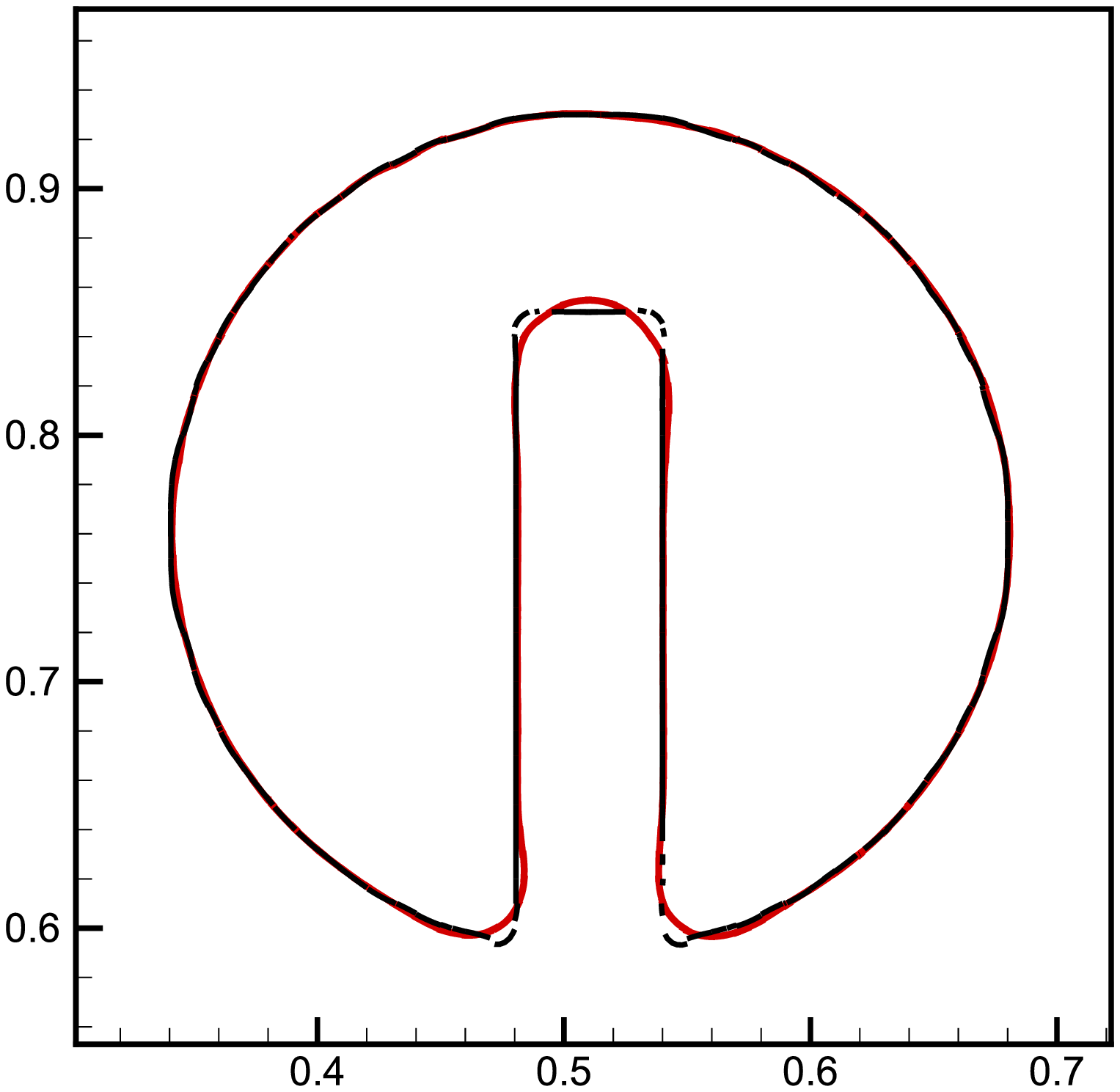}} \hspace{0.4cm}
	\subfigure[] {
		\centering
		\includegraphics[width=0.3\textwidth]{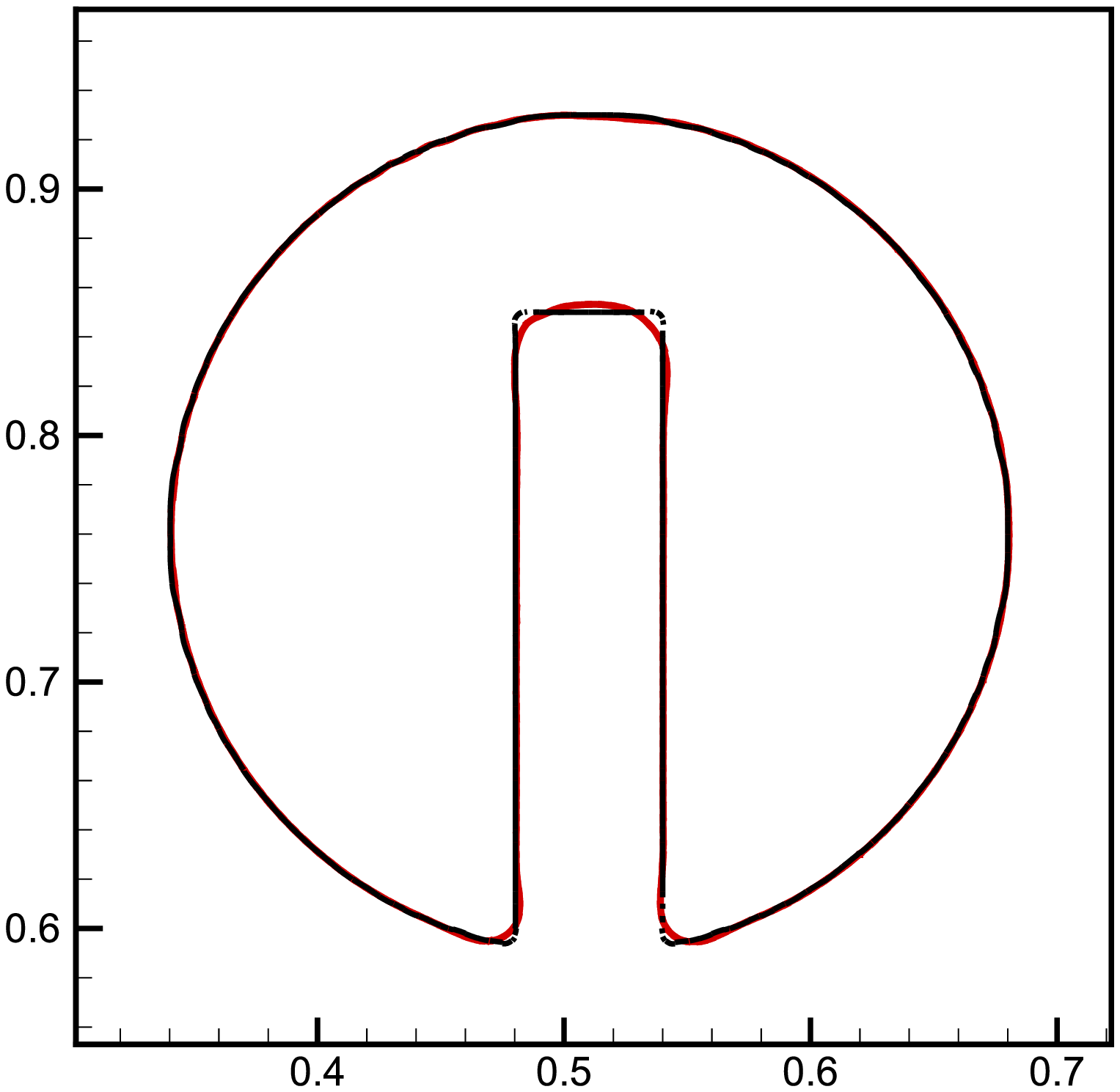}}
	\caption{The numerical results of Zalesak solid body rotation test after one revolution showing the reconstructed interfaces on (a) $50\times50$, (b) $100\times100$ and (c) $200\times200$ meshes. The PSI (red solid line) of the interface cells are plotted against the exact solution (black dashed line). }
\label{zalesak_reconstructed_results}
\end{figure}

We compute this test on structured grid for different grid sizes with 50, 100 and 200 vertices evenly distributed on each edge of computational domain. Fig.\ref{zalesak_results} shows the numerical results of interface identified by VOF 0.5 contour line.  As observed from these results, the interface in the slot region is well captured. As demonstrated in  \cite{xie2017toward} and \cite{qian2018}, the quadratic polynomial representation of the interface preserves the geometrical symmetry of the solution, which deteriorates significantly if a linear function (straight line) is used as commonly observed in the results of VOF method using PLIC reconstructions.  We also show the PSI of the interface cells in Fig. \ref{zalesak_reconstructed_results}. The interface is retrieved and represented by the cell-wise quadratic curves in the interface cells.    

\subsection{Vortex deformation transport test}

The THINC-scaling scheme is further assessed by single vortex test \cite{rider1998reconstructing}, in which a circle initially centred at $\left(0.5,0.75\right)$ in a unit square domain is advected by time dependent velocity field given by the stream function as follows,

\begin{equation}
\Psi\left(x,y,t\right)=\frac{1}{\pi}{\sin}^2\left(\pi{x}\right){\sin}^2\left(\pi{y}\right){\cos}\left(\frac{\pi{t}}{T}\right),
\end{equation}
where $T=8$ is specified in this test. This test, as one of the most widely used benchmark tests, is more challenging  to assess the capability of the scheme in capturing the heavily distorted interface with stretched tail when transported to $t=T/2$. From $t=T/2$ to $t=T$, the reverse velocity field restores the interface back to its initial shape. In case of $T=8$, the spiral tail becomes so thin that it can not be resolved by the resolution of a coarse grid. 
\begin{figure}[htbp]
	\centering
	\subfigure[] {
		\centering
		\includegraphics[width=0.3\textwidth]{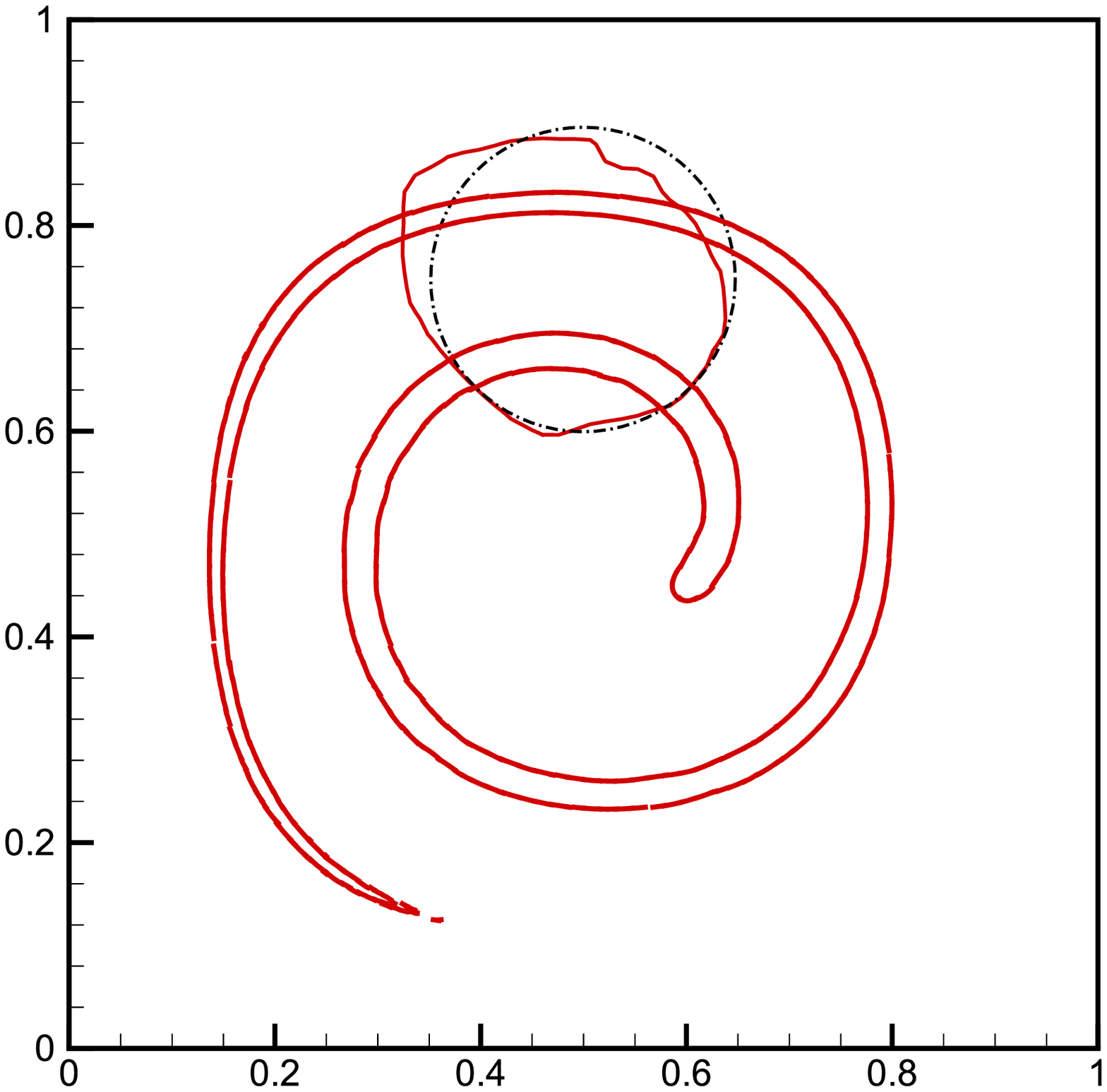}} \hspace{0.4cm}
	\subfigure[] {
		\centering
		\includegraphics[width=0.3\textwidth]{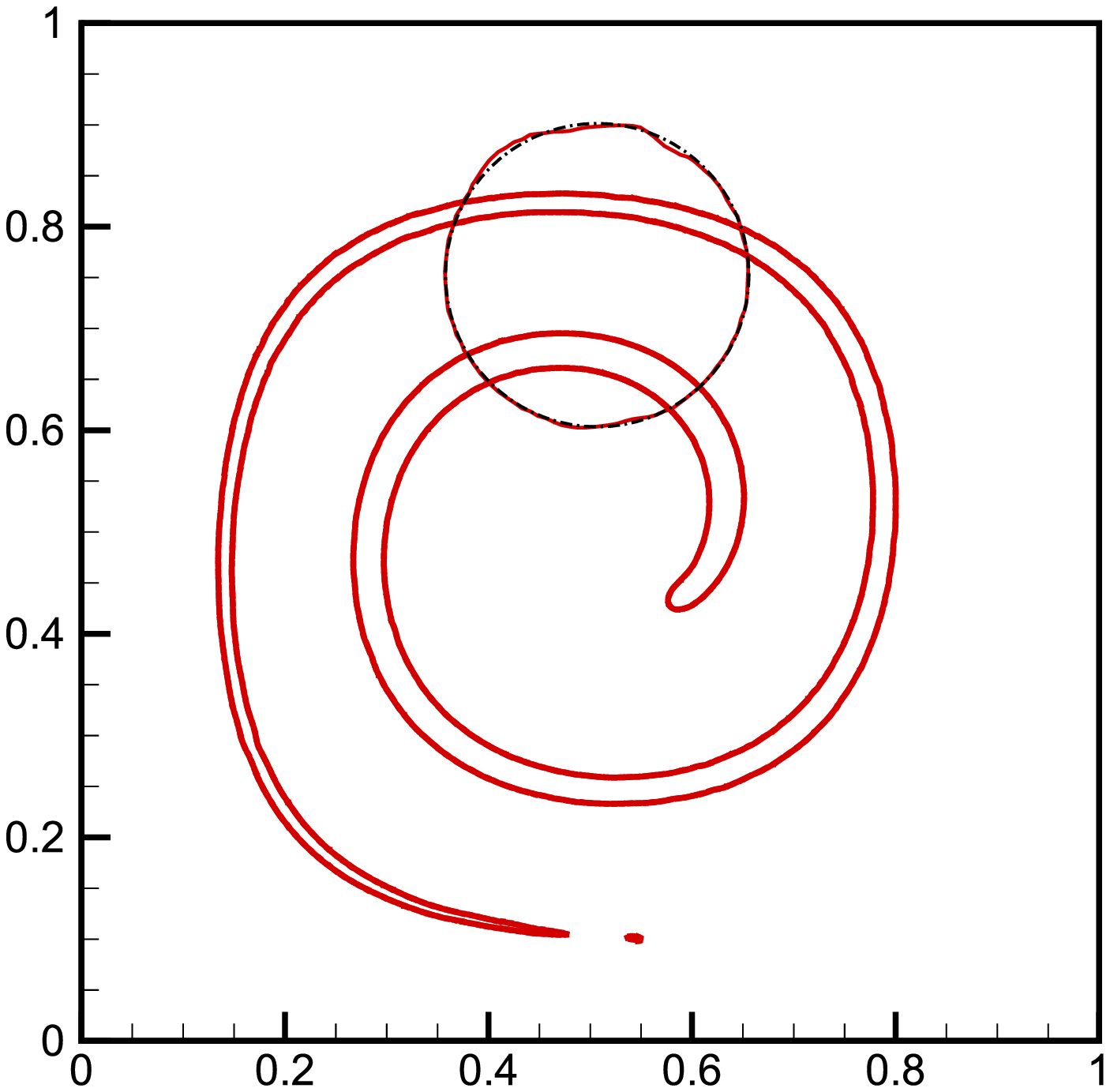}} \hspace{0.4cm}
	\subfigure[] {
		\centering
		\includegraphics[width=0.3\textwidth]{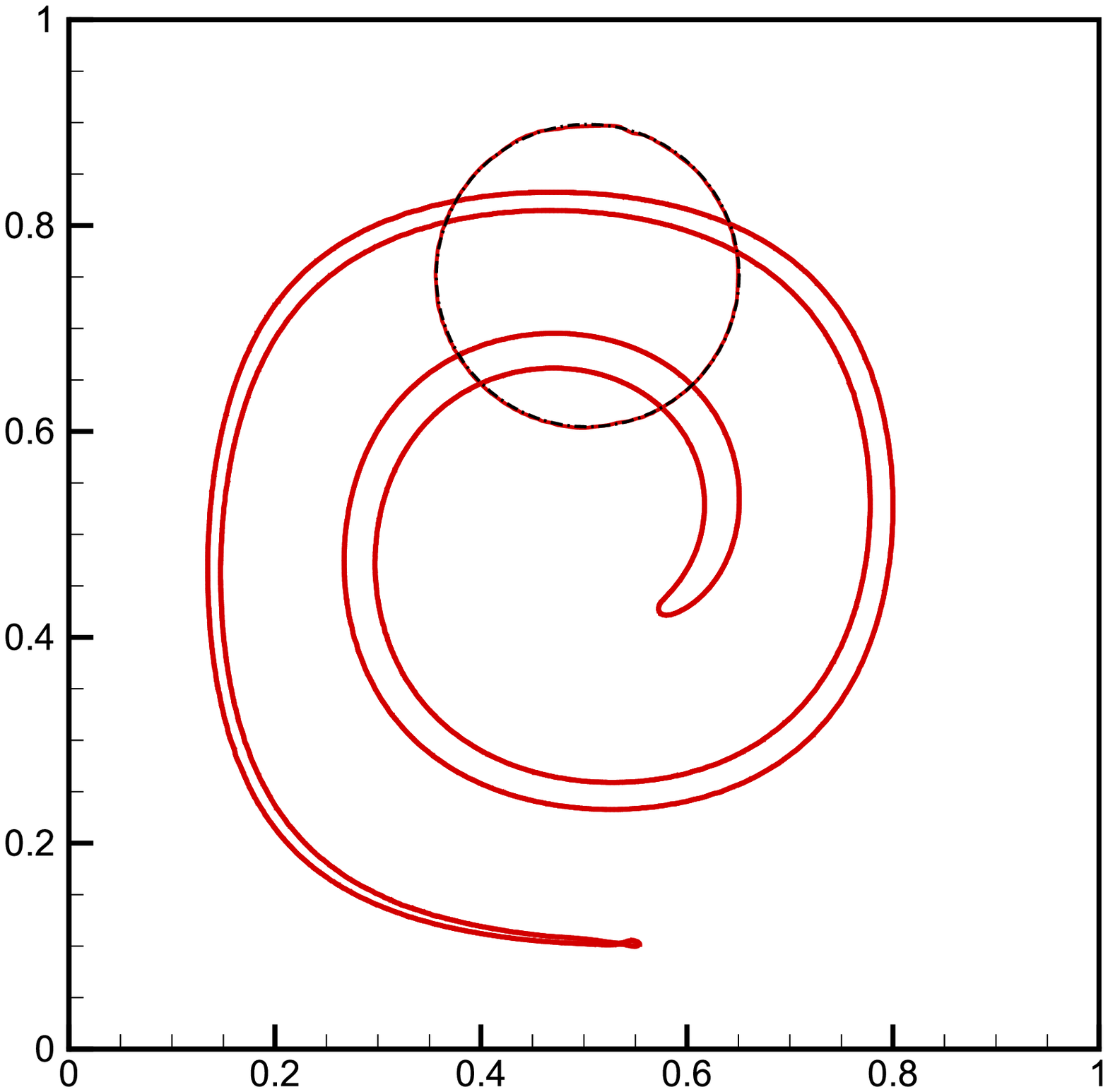}} \hspace{0.4cm}
	\caption{Numerical results for Rider-Kothe single vortex test showing the PSI for interface cells at $t=T/2$ and the VOF 0.5 contour lines at $t=T$ on (a) $64\times64$, (b) $128\times128$ and (c) $256\times256$ meshes}
\label{rider_results}
\end{figure}
\begin{figure}[htbp]
	\centering
	\subfigure[] {
		\centering
		\includegraphics[width=0.4\textwidth]{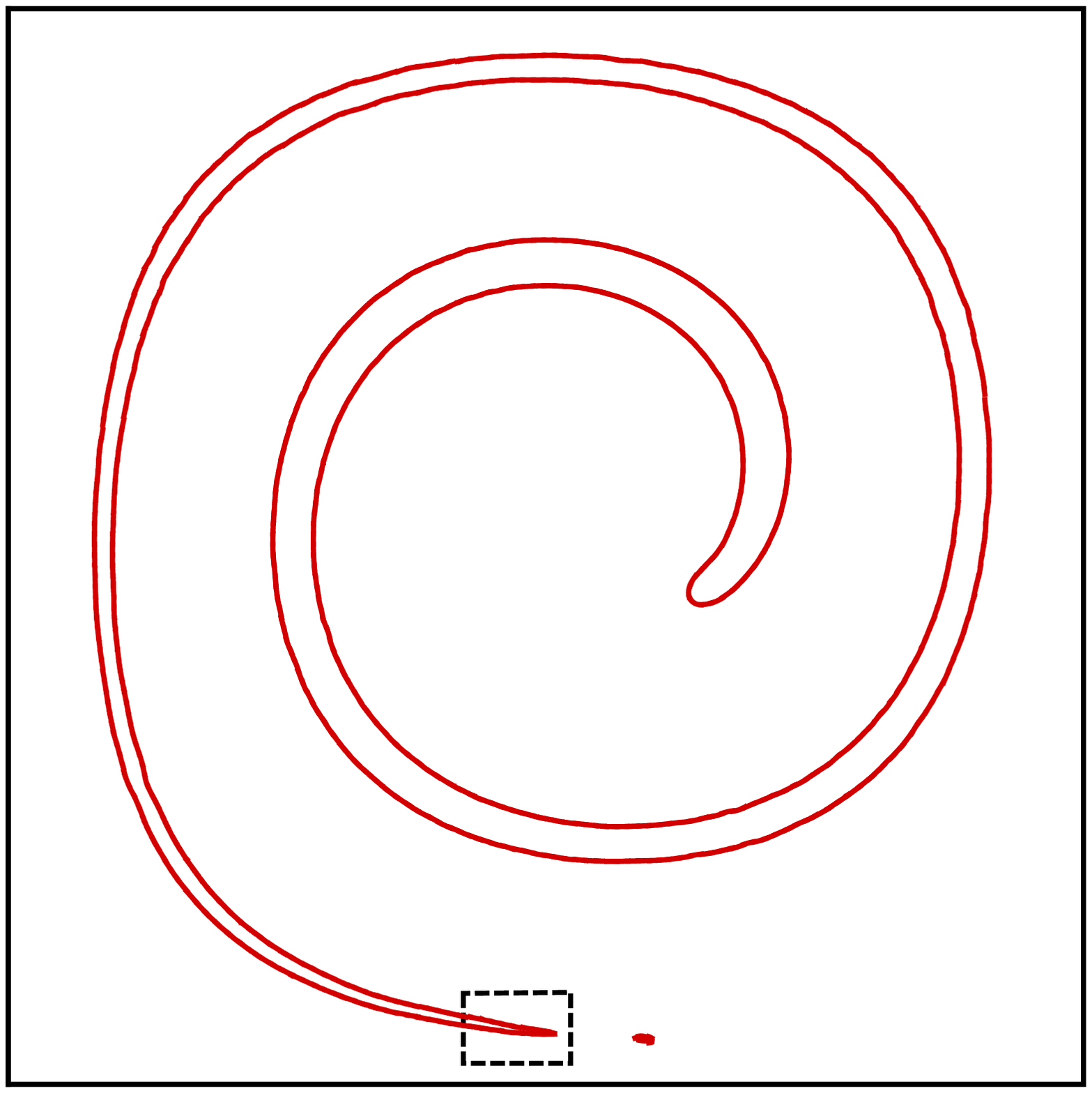}} \hspace{0.4cm}
	\centering
	\subfigure[] {
		\centering
		\includegraphics[width=0.4\textwidth]{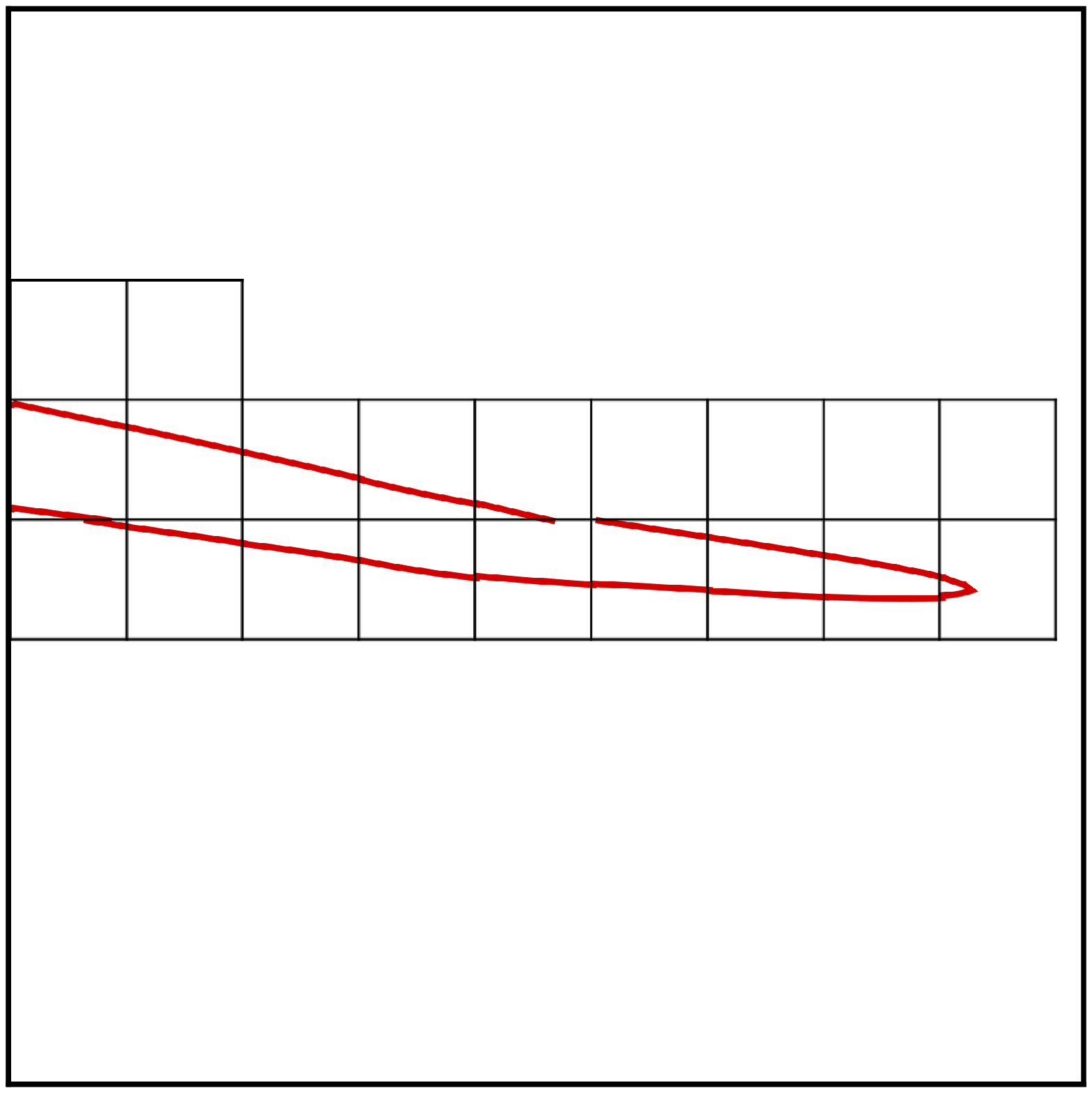}} \hspace{0.4cm}
	\centering
	\caption{The distorted interface at $T/2$ on $128\times128$ grid (a) where the boxed part is enlarged as the close-up of the thin tail in panel (b) which shows that the film of tail is under the resolution of grid cell, but can be still resolved by the THINC-scaling scheme. }
\label{thin-tail}
\end{figure}

\par We tested this case on Cartesian grids with different resolutions of $64\times64$,  $128\times128$ and $256\times256$ respectively. 
Numerical results on different grids at $t=T/2$ and $t=T$ are shown in Fig. \ref{rider_results}. The THINC-scaling scheme can capture the elongated tail even when the interface is under grid resolution, and can restore the initial circle with good solution quality. 

The PSI at $t=T/2$ on $128\times128$ grid is shown in Fig. \ref{thin-tail}. The tail tip is stretched into a thin film with a thickness smaller than the cell size. This sub-cell structure can still be reconstructed by the THINC-scaling scheme with quadratic or higher order polynomial representation. Consequently, the pieces of flotsam generated from the PLIC VOF methods are not observed here.

\begin{table}[ht]
\centering
\caption{Numerical errors and convergence rates for Rider-Kothe test on cartesian grid} 
\begin{tabular}[t]{lSSSSS} \toprule
{$\textbf{Methods}$}               & {$64^2$}              & {Order} & {$128^2$}               & {Order} & {$256^2$} \\ \midrule
{THINC/scaling}                    & {$1.28\times10^{-2}$} & {3.21}  & {$1.38\times10^{-3}$}   & {1.82}  & {$3.89\times10^{-4}$}\\
{UFVFC-Swartz \cite{Maric2018VOF}} & {$5.74\times10^{-3}$} & {1.98}  & {$1.45\times10^{-3}$}   & {1.94}  & {$3.77\times10^{-4}$}\\        
{Owkes and Desjardins \cite{Owkes2014VOF}}  & {$7.58\times10^{-3}$} & {2.01}  & {$1.88\times10^{-3}$}   & {2.21}  & {$4.04\times10^{-4}$}\\  
{isoAdvector-plicRDF \cite{Henning2019VOF}} & {$1.26\times10^{-2}$} & {2.27}  & {$2.61\times10^{-3}$}   & {2.19}  & {$5.71\times10^{-4}$}\\\bottomrule
\label{r-k-comparison}
\end{tabular}
\end{table}

A quantitative comparison with other sophisticated geometrical VOF methods \cite{Maric2018VOF,Owkes2014VOF,Henning2019VOF} is given in Table \ref{r-k-comparison}. It reveals the appealing accuracy of the present scheme.

 We solved the Rider-Kothe shear flow test case on an unstructured grid with triangular cells. In order to compare with the results in \cite{Henning2019VOF}, we set $64$ nodes on the domain boundaries in $x$ and $y$ directions respectively.  We used $\beta={6.0}/{\Delta}$, where $\Delta$ represents the cell size and is defined as the hydraulic diameter $\Delta={4A}/{P}$ with $A$ and $P$ being the area and perimeter of the triangular cell element respectively.

\begin{figure}[htbp]
	\centering
	\subfigure[]{
		\centering
		\includegraphics[width=0.3\textwidth]{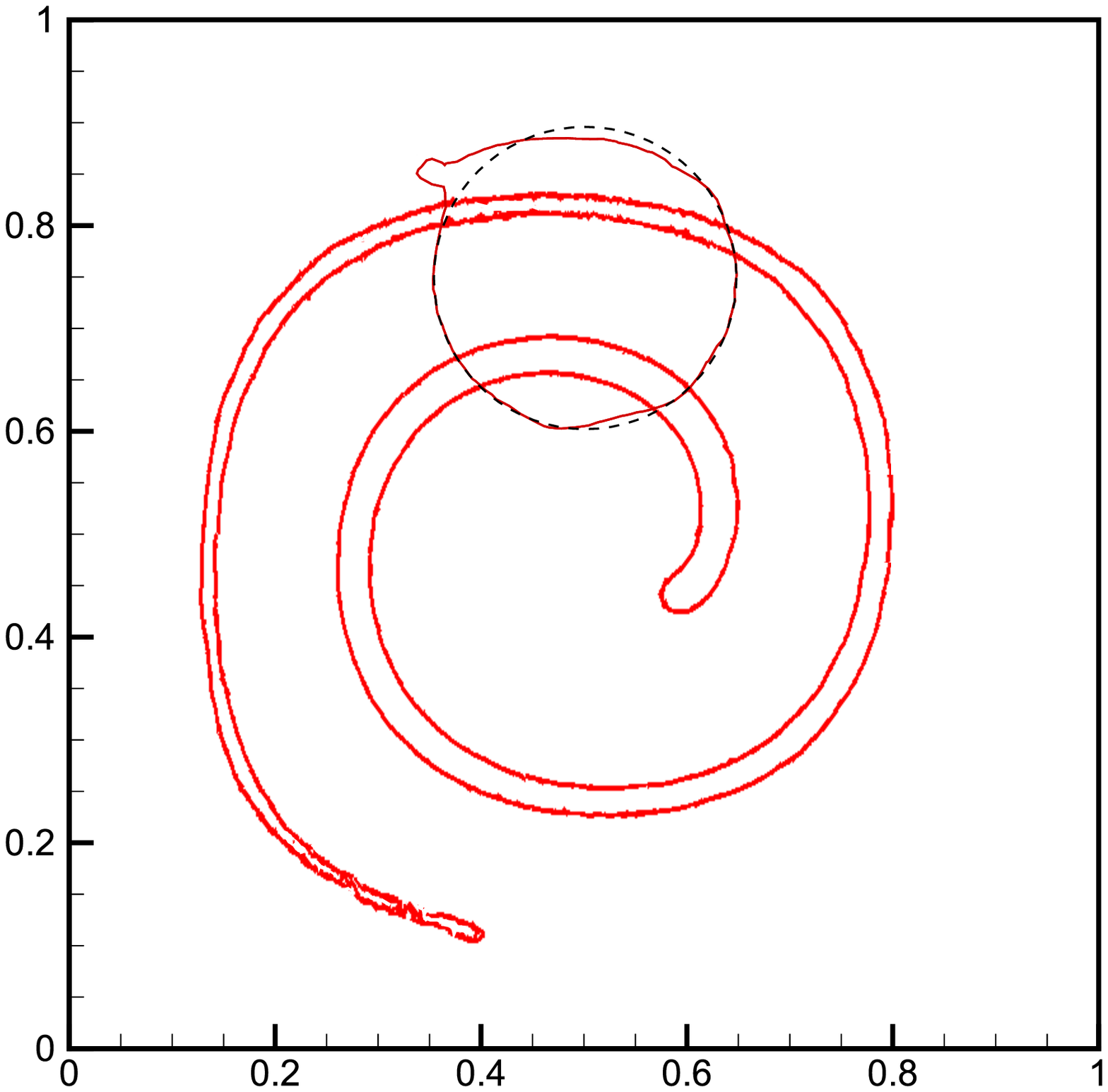}} \hspace{0.4cm}
	\subfigure[]{
		\centering
		\includegraphics[width=0.3\textwidth]{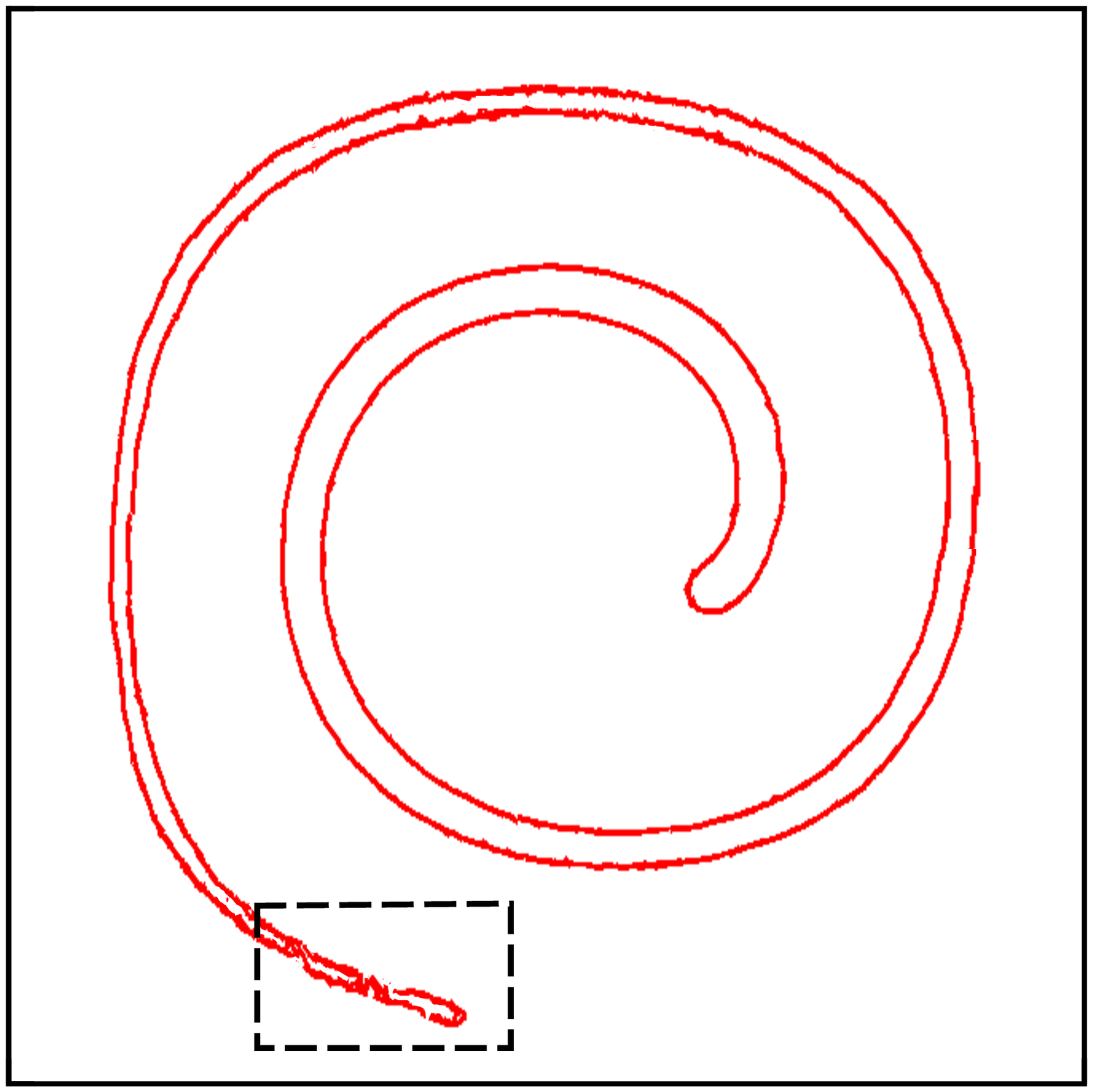}} \hspace{0.4cm}
			\subfigure[]{
		\centering
		\includegraphics[width=0.3\textwidth,height=0.4\textwidth]{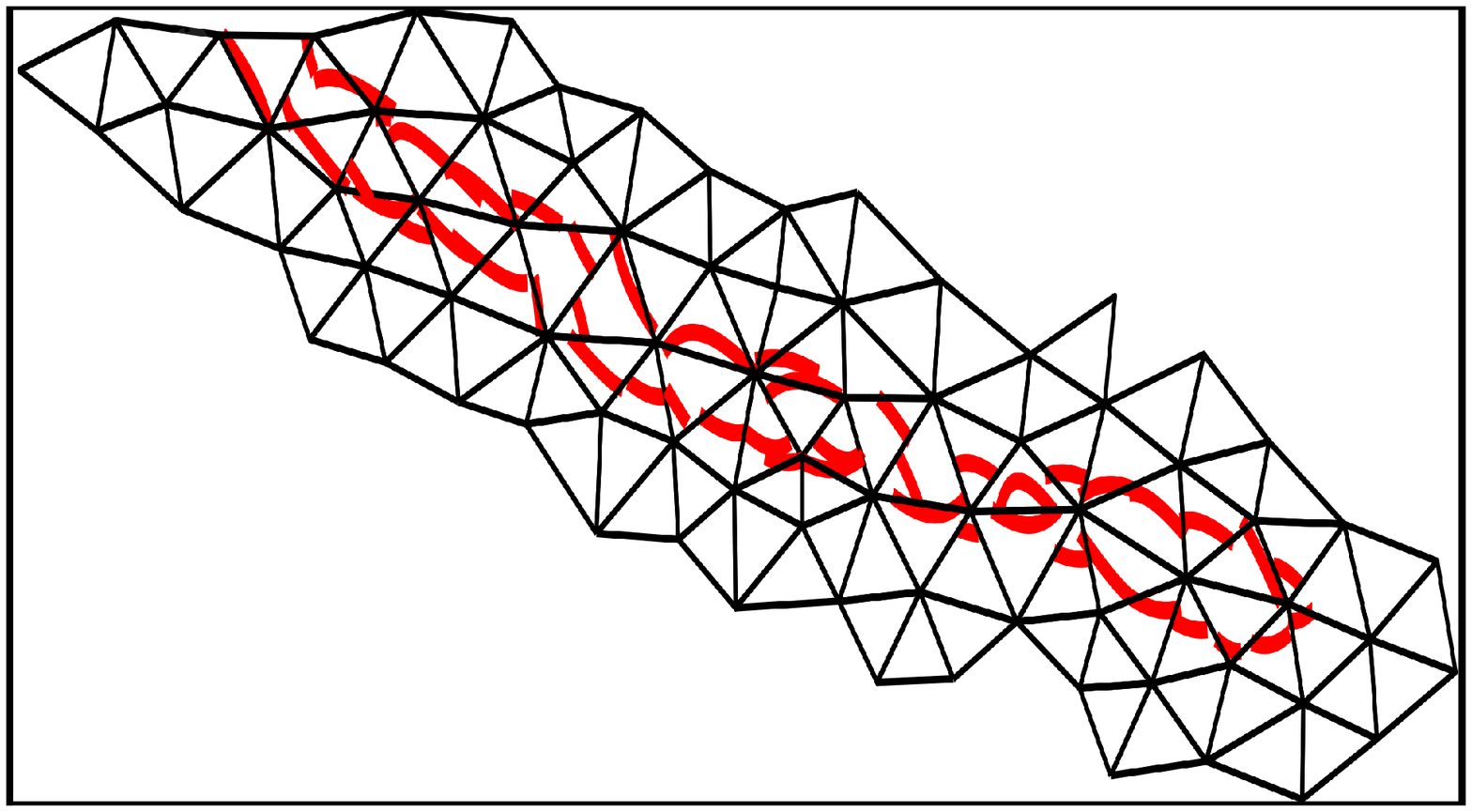}} \hspace{0.4cm}
	\caption{Numerical results for Rider-Kothe single vortex test on a triangular unstructured grid with $64$ nodes on each boundary of computational domain. (a): VOF 0.5 contour line at $t=T$ and  PSI of interface cells at $t=T/2$; (b): PSI of interface cells at $t=T/2$ with the stretched tail highlighted by the black dash-line box; (c): Enlarged view of PSI in the interface cells highlighted in (b). }
	\label{rider_unstructured}
\end{figure}

We plot the numerical results in  Fig. \ref{rider_unstructured}. The PSIs in the interface cells accurately present the reconstructed interface. The flotsams in the PLIC reconstruction \cite{Henning2019VOF} are avoided in the present results. A quantitative comparison is given in  Table \ref{r-k-unstructured-comparison}. The THINC-scaling scheme shows superiority in accuracy on unstructured grid. 

\begin{table}[ht]
\centering
\caption{Numerical errors and convergence rates for Rider-Kothe test on unstructured triangular grid with $64$ nodes on each boundary of computational domain.} 
\begin{tabular}[t]{lSSS} \toprule
{Methods}               & {THINC-scaling}                 & {THINC/QQ\cite{xie2017toward}}  & {isoAdvector-plicRDF \cite{Henning2019VOF}}\\ \midrule
{$L_1$ error}   & {$5.38\times10^{-3}$} & {$1.07\times10^{-2}$}     & {$2.21\times10^{-2}$}\\
\bottomrule
\label{r-k-unstructured-comparison}
\end{tabular}
\end{table}

\section{Conclusions}
We propose a novel scheme to unify the solution procedures of level set and VOF methods that are two interface-capturing methods based on completely different concept and numerical methodology. The underlying idea of the scheme, THINC-scaling scheme, is to use the THINC function to scale/convert between the level set function and a continuous Heaviside function which mimics the VOF field. 

The THINC function uses the level set polynomial to accurately retrieve the geometrical information of the interface from the level set field with high-order polynomials, and is constructed under the constraint of VOF value to rigorously satisfy the numerical conservativeness. Being a function handleable with conventional calculus tools, the THINC function facilitates efficient and accurate computations in interface-capturing schemes, which take the advantages from both VOF and level set methods. 

We verified the THINC-scaling scheme with advection benchmark tests for moving interfaces on both unstructured and unstructured grids, which demonstrate the super solution quality and the great potential of the proposed scheme as a moving interface-capturing scheme for practical utility.

\section{Acknowledgments}

This work was supported in part by the fund from JSPS (Japan Society for the Promotion of Science) under Grant No. 18H01366. 

\clearpage{}
\bibliographystyle{elsarticle-num-names}
\bibliography{THINC_LS} 

\begin{thebibliography}{25}
\providecommand{\natexlab}[1]{#1}
\providecommand{\url}[1]{\texttt{#1}}
\providecommand{\urlprefix}{URL }
\expandafter\ifx\csname urlstyle\endcsname\relax
  \providecommand{\doi}[1]{doi:\discretionary{}{}{}#1}\else
  \providecommand{\doi}[1]{doi:\discretionary{}{}{}\begingroup
  \urlstyle{rm}\url{#1}\endgroup}\fi
\providecommand{\bibinfo}[2]{#2}

\bibitem[{Hirt and Nichols(1981)}]{hirt1981volume}
\bibinfo{author}{C.~W. Hirt}, \bibinfo{author}{B.~D. Nichols},
  \bibinfo{title}{Volume of fluid (VOF) method for the dynamics of free
  boundaries}, \bibinfo{journal}{Journal of computational physics}
  \bibinfo{volume}{39}~(\bibinfo{number}{1}) (\bibinfo{year}{1981})
  \bibinfo{pages}{201--225}.

\bibitem[{Youngs(1982)}]{youngs1982time}
\bibinfo{author}{D.~L. Youngs}, \bibinfo{title}{Time-dependent multi-material
  flow with large fluid distortion}, \bibinfo{journal}{Numerical methods for
  fluid dynamics} .

\bibitem[{Lafaurie et~al.(1994)Lafaurie, Nardone, Scardovelli, Zaleski, and
  Zanetti}]{lafaurie1994modelling}
\bibinfo{author}{B.~Lafaurie}, \bibinfo{author}{C.~Nardone},
  \bibinfo{author}{R.~Scardovelli}, \bibinfo{author}{S.~Zaleski},
  \bibinfo{author}{G.~Zanetti}, \bibinfo{title}{Modelling merging and
  fragmentation in multiphase flows with SURFER}, \bibinfo{journal}{Journal of
  Computational Physics} \bibinfo{volume}{113}~(\bibinfo{number}{1})
  (\bibinfo{year}{1994}) \bibinfo{pages}{134--147}.

\bibitem[{Rider and Kothe(1998)}]{rider1998reconstructing}
\bibinfo{author}{W.~J. Rider}, \bibinfo{author}{D.~B. Kothe},
  \bibinfo{title}{Reconstructing volume tracking}, \bibinfo{journal}{Journal of
  computational physics} \bibinfo{volume}{141}~(\bibinfo{number}{2})
  (\bibinfo{year}{1998}) \bibinfo{pages}{112--152}.

\bibitem[{Scardovelli and Zaleski(2000)}]{scardovelli2000analytical}
\bibinfo{author}{R.~Scardovelli}, \bibinfo{author}{S.~Zaleski},
  \bibinfo{title}{Analytical relations connecting linear interfaces and volume
  fractions in rectangular grids}, \bibinfo{journal}{Journal of Computational
  Physics} \bibinfo{volume}{164}~(\bibinfo{number}{1}) (\bibinfo{year}{2000})
  \bibinfo{pages}{228--237}.

\bibitem[{Osher and Sethian(1988)}]{osher1988fronts}
\bibinfo{author}{S.~Osher}, \bibinfo{author}{J.~A. Sethian},
  \bibinfo{title}{Fronts propagating with curvature-dependent speed: algorithms
  based on Hamilton-Jacobi formulations}, \bibinfo{journal}{Journal of
  computational physics} \bibinfo{volume}{79}~(\bibinfo{number}{1})
  (\bibinfo{year}{1988}) \bibinfo{pages}{12--49}.

\bibitem[{Sethian(1999)}]{sethian1999level}
\bibinfo{author}{J.~A. Sethian}, \bibinfo{title}{Level set methods and fast
  marching methods: evolving interfaces in computational geometry, fluid
  mechanics, computer vision, and materials science}, vol.~\bibinfo{volume}{3},
  \bibinfo{publisher}{Cambridge university press}, \bibinfo{year}{1999}.

\bibitem[{Osher and Fedkiw(2003)}]{osher2003implicit}
\bibinfo{author}{S.~Osher}, \bibinfo{author}{R.~Fedkiw},
  \bibinfo{title}{Implicit Functions}, in: \bibinfo{booktitle}{Level Set
  Methods and Dynamic Implicit Surfaces}, \bibinfo{publisher}{Springer},
  \bibinfo{pages}{3--16}, \bibinfo{year}{2003}.

\bibitem[{Sussman and Puckett(2000)}]{sussman2000coupled}
\bibinfo{author}{M.~Sussman}, \bibinfo{author}{E.~G. Puckett},
  \bibinfo{title}{A coupled level set and volume-of-fluid method for computing
  3D and axisymmetric incompressible two-phase flows},
  \bibinfo{journal}{Journal of computational physics}
  \bibinfo{volume}{162}~(\bibinfo{number}{2}) (\bibinfo{year}{2000})
  \bibinfo{pages}{301--337}.

\bibitem[{M{\'e}nard et~al.(2007)M{\'e}nard, Tanguy, and
  Berlemont}]{menard2007coupling}
\bibinfo{author}{T.~M{\'e}nard}, \bibinfo{author}{S.~Tanguy},
  \bibinfo{author}{A.~Berlemont}, \bibinfo{title}{Coupling level set/VOF/ghost
  fluid methods: Validation and application to 3D simulation of the primary
  break-up of a liquid jet}, \bibinfo{journal}{International Journal of
  Multiphase Flow} \bibinfo{volume}{33}~(\bibinfo{number}{5})
  (\bibinfo{year}{2007}) \bibinfo{pages}{510--524}.

\bibitem[{Yang et~al.(2006)Yang, James, Lowengrub, Zheng, and
  Cristini}]{yang2006adaptive}
\bibinfo{author}{X.~Yang}, \bibinfo{author}{A.~J. James},
  \bibinfo{author}{J.~Lowengrub}, \bibinfo{author}{X.~Zheng},
  \bibinfo{author}{V.~Cristini}, \bibinfo{title}{An adaptive coupled
  level-set/volume-of-fluid interface capturing method for unstructured
  triangular grids}, \bibinfo{journal}{Journal of Computational Physics}
  \bibinfo{volume}{217}~(\bibinfo{number}{2}) (\bibinfo{year}{2006})
  \bibinfo{pages}{364--394}.

\bibitem[{Sun and Tao(2010)}]{sun2010coupled}
\bibinfo{author}{D.~Sun}, \bibinfo{author}{W.~Tao}, \bibinfo{title}{A coupled
  volume-of-fluid and level set (VOSET) method for computing incompressible
  two-phase flows}, \bibinfo{journal}{International Journal of Heat and Mass
  Transfer} \bibinfo{volume}{53}~(\bibinfo{number}{4}) (\bibinfo{year}{2010})
  \bibinfo{pages}{645--655}.

\bibitem[{Xiao et~al.(2005)Xiao, Honma, and Kono}]{xiao2005simple}
\bibinfo{author}{F.~Xiao}, \bibinfo{author}{Y.~Honma},
  \bibinfo{author}{T.~Kono}, \bibinfo{title}{A simple algebraic interface
  capturing scheme using hyperbolic tangent function},
  \bibinfo{journal}{International Journal for Numerical Methods in Fluids}
  \bibinfo{volume}{48}~(\bibinfo{number}{9}) (\bibinfo{year}{2005})
  \bibinfo{pages}{1023--1040}.

\bibitem[{Xiao et~al.(2011)Xiao, Ii, and Chen}]{xiao2011revisit}
\bibinfo{author}{F.~Xiao}, \bibinfo{author}{S.~Ii}, \bibinfo{author}{C.~Chen},
  \bibinfo{title}{Revisit to the THINC scheme: a simple algebraic VOF
  algorithm}, \bibinfo{journal}{Journal of Computational Physics}
  \bibinfo{volume}{230}~(\bibinfo{number}{19}) (\bibinfo{year}{2011})
  \bibinfo{pages}{7086--7092}.

\bibitem[{Ii et~al.(2012)Ii, Sugiyama, Takeuchi, Takagi, Matsumoto, and
  Xiao}]{ii2012interface}
\bibinfo{author}{S.~Ii}, \bibinfo{author}{K.~Sugiyama},
  \bibinfo{author}{S.~Takeuchi}, \bibinfo{author}{S.~Takagi},
  \bibinfo{author}{Y.~Matsumoto}, \bibinfo{author}{F.~Xiao}, \bibinfo{title}{An
  interface capturing method with a continuous function: The THINC method with
  multi-dimensional reconstruction}, \bibinfo{journal}{Journal of Computational
  Physics} \bibinfo{volume}{231}~(\bibinfo{number}{5}) (\bibinfo{year}{2012})
  \bibinfo{pages}{2328--2358}.

\bibitem[{Xie and Xiao(2017)}]{xie2017toward}
\bibinfo{author}{B.~Xie}, \bibinfo{author}{F.~Xiao}, \bibinfo{title}{Toward
  efficient and accurate interface capturing on arbitrary hybrid unstructured
  grids: The THINC method with quadratic surface representation and Gaussian
  quadrature}, \bibinfo{journal}{Journal of Computational Physics}
  \bibinfo{volume}{349} (\bibinfo{year}{2017}) \bibinfo{pages}{415--440}.

\bibitem[{Qian et~al.(2018)Qian, Wei, and Xiao}]{qian2018}
\bibinfo{author}{L.~Qian}, \bibinfo{author}{Y.~Wei}, \bibinfo{author}{F.~Xiao},
  \bibinfo{title}{Coupled THINC and level set method: A conservative interface
  capturing scheme with high-order surface representations},
  \bibinfo{journal}{Journal of Computational Physics} \bibinfo{volume}{373}
  (\bibinfo{year}{2018}) \bibinfo{pages}{284--303}.

\bibitem[{Shu(1988)}]{shu88}
\bibinfo{author}{C.-W. Shu}, \bibinfo{title}{Total-variation-diminishing time
  discretizations}, \bibinfo{journal}{SIAM J. Sci. Stat. Comput.}
  \bibinfo{volume}{9} (\bibinfo{year}{1988}) \bibinfo{pages}{1073--1084}.

\bibitem[{Zhao(2005)}]{zhao2005fast}
\bibinfo{author}{H.~Zhao}, \bibinfo{title}{A fast sweeping method for eikonal
  equations}, \bibinfo{journal}{Mathematics of computation}
  \bibinfo{volume}{74}~(\bibinfo{number}{250}) (\bibinfo{year}{2005})
  \bibinfo{pages}{603--627}.

\bibitem[{Dianat et~al.(2017)Dianat, Skarysz, and
  Garmory}]{Dianat2017unstructured}
\bibinfo{author}{M.~Dianat}, \bibinfo{author}{M.~Skarysz},
  \bibinfo{author}{A.~Garmory}, \bibinfo{title}{A Coupled Level Set and Volume
  of Fluid method for automotive exterior water management applications},
  \bibinfo{journal}{International Journal of Multiphase Flow}
  \bibinfo{volume}{91} (\bibinfo{year}{2017}) \bibinfo{pages}{19--38}.

\bibitem[{Strain(1999)}]{strain99}
\bibinfo{author}{J.~Strain}, \bibinfo{title}{Semi-Lagrange methods for level
  set equations}, \bibinfo{journal}{Journal of Computational Physics}
  \bibinfo{volume}{2} (\bibinfo{year}{1999}) \bibinfo{pages}{498--533}.

\bibitem[{Zalesak(1979)}]{zalesak1979fully}
\bibinfo{author}{S.~T. Zalesak}, \bibinfo{title}{Fully multidimensional
  flux-corrected transport algorithms for fluids}, \bibinfo{journal}{Journal of
  computational physics} \bibinfo{volume}{31}~(\bibinfo{number}{3})
  (\bibinfo{year}{1979}) \bibinfo{pages}{335--362}.

\bibitem[{Maric et~al.(2018)Maric, Marschall, and Bothe}]{Maric2018VOF}
\bibinfo{author}{T.~Maric}, \bibinfo{author}{H.~Marschall},
  \bibinfo{author}{D.~Bothe}, \bibinfo{title}{An enhanced un-split face-vertex
  flux-based VoF method}, \bibinfo{journal}{Journal of computational physics}
  \bibinfo{volume}{371} (\bibinfo{year}{2018}) \bibinfo{pages}{967--993}.

\bibitem[{Owkes and Desjardins(2014)}]{Owkes2014VOF}
\bibinfo{author}{M.~Owkes}, \bibinfo{author}{O.~Desjardins}, \bibinfo{title}{A
  computational framework for conservative, three-dimensional, unsplit,
  geometric transport with application to the volume-of-fluid (VOF) method},
  \bibinfo{journal}{Journal of computational physics} \bibinfo{volume}{270}
  (\bibinfo{year}{2014}) \bibinfo{pages}{587--612}.

\bibitem[{Scheufler and Roenby(2018)}]{Henning2019VOF}
\bibinfo{author}{H.~Scheufler}, \bibinfo{author}{J.~Roenby},
  \bibinfo{title}{Accurate and efficient surface reconstruction from volume
  fraction data on general meshes}, \bibinfo{journal}{Journal of computational
  physics} \bibinfo{volume}{383} (\bibinfo{year}{2018}) \bibinfo{pages}{1--23}.

\end{thebibliography}

\end{document}